\newcommand{\free}{\text{\huge{$\ast$}}}

\newcommand{\mul}{Mul[[\mathfrak{B}]]}

\documentclass[10pt,draft,reqno]{amsart}
    \makeatletter
     \def\section{\@startsection{section}{1}%
     \z@{.7\linespacing\@plus\linespacing}{.5\linespacing}%
     {\bfseries%\normalfont\scshape
     \centering
     }}
     \def\@secnumfont{\bfseries}
     \makeatother
\newtheorem{theorem}{Theorem}[section]
\newtheorem{lemma}[theorem]{Lemma}

\newtheorem{corollary}[theorem]{Corollary}
\theoremstyle{definition}
\newtheorem{definition}[theorem]{Definition}

\theoremstyle{remark}
\newtheorem{remark}[theorem]{Remark}
\numberwithin{equation}{section} 
\usepackage{amssymb, graphics}
\begin{document}

\title[Multilinear function series, c-free probability with
amalgamation] {Multilinear function series in  conditionally free
probability with amalgamation}

\author{Mihai Popa}
\address{Mihai Popa: Department of Mathematics,
Indiana University at Bloomington, \linebreak {Rawles Hall},
 931 E 3rd St, Bloomington, IN 47405}
\email{mipopa@indiana.edu}

\subjclass[2000] {Primary 45L53; Secondary  46L08}

\keywords{ conditional freeness, conditional expectation,
R-transform,
 multilinear function series}

\begin{abstract}
As in the cases of freeness and monotonic independence, the notion
of
 conditional freeness is meaningful when complex-valued states are
 replaced by positive conditional expectations. In this framework, the paper presents
 several positivity results, a version of the central limit theorem
 and an analogue of the conditionally free R-transform constructed by
 means of multilinear function series.
\end{abstract}

\maketitle

\section{Introduction}

The paper addresses a topic related to conditionally free (or,
shortly, using the term from \cite{bls}, c-free) probability. This
notion was developed in the '90's  (see \cite{bs}, \cite{bls}) as an
extension of freeness within the framework of $\ast$-algebras
endowed with not one, but two states. Namely, given a family of
unital algebras $\{\mathfrak{A}\}_{i\in \mathcal{I}}$, each
$\mathfrak{A}_i$ endowed with two expectations $\varphi_i,
\psi_i:\mathfrak{A}_i\longrightarrow\mathbb{C}$, their c-free
product is the triple $(\mathfrak{A}, \varphi, \psi)$, where:
\begin{enumerate}
\item[(i)] $\mathfrak{A}=\free_{i\in I}\mathfrak{A}_i$ is the free
product of the algebras $\mathfrak{A}_i$.
\item[(ii)] $\psi=\free_{i\in \mathcal{I}}\psi_i$ and
$\varphi=\free_{(\psi_i)_{i\in \mathcal{I}}}\varphi_i$ are
expectations given by the relations
\begin{enumerate}
\item[(a)]$\psi(a_1\cdots a_n)=0$
\item[(b)]$\varphi(a_1\cdots a_n)=\varphi_{\varepsilon(1)}(a_1)
\cdots\varphi_{\varepsilon(n)}(a_n)$
\end{enumerate}
\end{enumerate}
for all $a_j\in\mathfrak{A}_{\varepsilon(j)}, j=1,\dots,n$ such that
$\psi_{\varepsilon(j)}(a_j)=0$ and $\varepsilon(1)\neq\dots\neq
\varepsilon(n)$.

A key result is that if the $\mathfrak{A}_i$ are $\ast$-algebras and
$\varphi_i, \psi_i$ are positive functionals, then $\varphi$ and
$\psi$ are also positive.

In \cite{mlotk1}, the positivity of the free product maps $\varphi,
\psi$ is proved for the case when $\varphi_1,\varphi_2$ are positive
conditional expectations in a common $C*$-subalgebra, but
$\psi_1,\psi_2$ remain positive $\mathbb{C}$-valued maps. A more
general situation is indeed discussed (see Theorem 3, Section 6,
from \cite{mlotk1}), but the question if $\varphi, \psi$ re positive
for $\varphi_{1,2}, \psi_{1,2}$ arbitrary positive conditional
expectations is left unanswered.

 A first answer was given in \cite{mvp}, where we showed
  that for $\mathfrak{A}$ a $\ast$-algebra, the
analogous construction with both $\varphi$ and $\psi$ valued in a
$C^*$-subalgebra $\mathfrak{B}$ of $\mathfrak{A}$ still retains the
positivity. The present paper further develops  this result (see
Theorem \ref{condpositivity})  and also demonstrates the use of
multilinear function series in c-free setting.

\pagebreak

In \cite{bls} is constructed a c-free version of Voiculescu's
$R$-transform, which we will call the ${}^cR$-transform, with the
property that ${}^cR_{X+Y}={}^cR_{X}+{}^cR_{Y}$ if $X$ and $Y$ are
c-free elements from the algebra $\mathfrak{A}$ relative to
$\varphi$ and $\psi$ (i.e. the relations (a) and (b) from the
definition of the c-free product hold true for the subalgebras
generated by $X$ and $Y$.)

The apparatus of multilinear function series is used in  recent work
of K. Dykema (\cite{dykema} and \cite{dykema2}) to construct
suitable analogues for the $R$ and $S$-transforms in the framework
of freeness with amalgamation. We will show that this construction
is also appropriate for the ${}^cR$-transform mentioned above. The
techniques used differ from the ones of \cite{dykema}, the Fock
space type construction being substituted by combinatorial
techniques similar to \cite{bls} and \cite{ns}. Particularly,
Theorems 3.3 and 3.6 contain new (shorter) proves of the results
6.1--6.13 from \cite{dykema}.

The paper is structured in four sections. In Section 2 are stated
the basic definitions and  are proved the main positivity results.
Section 3 describes the construction and the basic property of the
multilinear function series ${}^cR$-transform and Section 4 treats
the central limit theorem and a related positivity property.
%
%%%%%%%%%%%%%%%%%%%%%%%%%%%%%%%%%%%%%%%%%%%%%%%%%%%%%%%%%%%%%%%%%%%%%%%%%%%%%%%%%%%%%%%
%2222222222222222222222222222222222222222222222222222222222222222222222222
%
%
\section{Definitions and positivity results}\label{second}
\begin{definition}\label{defncfree}
Let $\{\mathfrak{A}_i\}_{i\in{\mathcal{I}}}$ be a family of
algebras, all containing the subalgebra $\mathfrak{B}$. Suppose
$\mathfrak{D}$ is a subalgebra of $\mathfrak{B}$ and
$\Psi_i:\mathfrak{A}_i\longrightarrow\mathfrak{D}$ and
$\Phi_i:\mathfrak{A}_i\longrightarrow\mathfrak{B}$ are conditional
expectations, $i \in\mathcal{I}$. We say that the triple
$(\mathfrak{A},
\Phi,\Psi)=\free_{i\in{I}}(\mathfrak{A}_i,\Phi_i,\Psi_i)$ is the
{conditionally free product} with amalgamation over
($\mathfrak{B}$,$\mathfrak{D}$), or shortly, the \emph{c-free
product}, of the triples
$(\mathfrak{A}_i,\Phi_i,\Psi_i)_{{i\in{I}}}$ if
\begin{enumerate}
\item[(1)]{$\mathfrak{A}$ is the free product
with amalgamation over $\mathfrak{B}$ of the family
$(\mathfrak{A}_i)_{i\in\mathcal{I}}$}
\item[(2)]{$\Psi=\free_{i\in\mathcal{I}}\Psi_i$ and
$\Phi=\free_{(\Psi_i)_{i\in\mathcal{I}}}\Phi_i$ are determined by
the relations}
\begin{eqnarray*}
\Psi(a_1a_2\dots a_n)&=&0\\
\Phi(a_1a_2\dots a_n)&=&\Phi(a_1)\Phi(a_2)\dots \Phi(a_n),
\end{eqnarray*}
\end{enumerate}
{for all $a_i\in\mathfrak{A}_{\varepsilon(i)},
\varepsilon(i)\in\mathcal{I}$\ such that
$\varepsilon(1)\neq\varepsilon(2)\neq\dots\neq \varepsilon(n)$\ and
$\Psi_{\varepsilon(i)}(a_i)=0$.}

\end{definition}

When $\mathfrak{D}=\mathbb{C}$, this definition reduces to the one
given in \cite{mlotk1}. When both $\mathfrak{B}$ and $\mathfrak{D}$
are equal to $\mathbb{C}$, this definition was given in \cite{bls}.

 When discussing positivity, we need a $\ast$-structure on our algebras.
  We will demand that $\mathfrak{B}$
 and $\mathfrak{D}$ be C$^*$-algebras, while $\mathfrak{A}_i$ and $\mathfrak{A}$ are only
 required to be $\ast$-algebras.

 The following results are slightly modified versions of Lemma 6.4
 and Theorem 6.5 from \cite{mvp}.

 \begin{lemma}\label{mixprod}
 Let $\mathfrak{B}$ be a  $C^*$-algebra and $\mathfrak{A}_1$, $\mathfrak{A}_2$  be two $\ast$-algebras
 containing $\mathfrak{B}$ as a $\ast$-subalgebra, endowed with  positive conditional
 expectations
 $\Phi_j:\mathfrak{A}_j\longrightarrow\mathfrak{B}$, $j=1,2$.
 If $a_1,\dots,a_n\in\mathfrak{A}_1$, $a_{n+1},\dots,a_{n+m}\in\mathfrak{A}_2$
 and
  $A=(A_{i,j})_{i,j}\in M_{n+m}(\mathfrak{B})$
  is the matrix with the entries
 \[
  A_{i,j}=\left\{
 \begin{matrix}
 \Phi_1(a_i^*a_j)& \text{if}& i,j\leq n
  \hspace{1.2cm} \\
 \Phi_1(a_i^*)\Phi_2(a_j)& \text{if}& i\leq n, j>n\\
\Phi_2(a_i^*)\Phi_1(a_j)& \text{if}& i>n, j\leq n\\
\Phi_2(a_i^*a_j)& \text{if}& i,j> n\hspace{1.2cm}\\
 \end{matrix}\right
  .\]
 then $A$ is positive.
 \end{lemma}
 \begin{proof}
 The vector space
 $\mathfrak{E}=\mathfrak{B}\oplus \ker(\Phi_1)\oplus \ker(\Phi_2)$
 has a $\mathfrak{B}$-bimodule structure given by the algebraic operations on
 $\mathfrak{A}_1$ and $\mathfrak{A}_2$. Consider the $\mathfrak{B}$-sesquilinear pairing
 \[
  \langle\cdot,\cdot\rangle:\mathfrak{E}\times\mathfrak{E}\longrightarrow\mathfrak{B}
  \]
determined by the relations:
 \begin{eqnarray*}
 \langle b_1, b_2\rangle&=&b_1^*b_2,\ \ \text{for}\
  b_1,
 b_2\in\mathfrak{B}\\
 \langle u_j,v_j\rangle&=&\Phi_j(u_j^*v_j),\ \ \text{for}\
  u_j,
 v_j\in\ker(\Phi_j), j=1,2\\
 \langle u_1,u_2\rangle&=&\langle u_2,u_1\rangle=0\ \ \text{for}\
  u_1\in \ker(\Phi_1),\ \text{and}\
 u_2\in\ker(\Phi_2).\\
 \langle b,u_j\rangle&=&\langle u_j, b\rangle=0\ \
 \text{for all}\ b\in \mathfrak{B},u_j\in\mathfrak{A}_j
 \end{eqnarray*}

With this notation, we have that
 $A_{i,j}=\langle a_i,a_j\rangle$,
 hence it
suffices to show that  $\langle a,a\rangle\geq 0$ for all
$a\in\mathfrak{E}$.

Indeed, for an element $a=b+u_1+u_2$ with $b\in\mathfrak{B},
u_j\in\ker(\Phi_j), j=1,2$, we have:
\begin{eqnarray*}
\langle a,a\rangle &=& \langle b+u_1+u_2, b+u_1+u_2\rangle\\
&=&\langle b, b\rangle+\langle u_1,u_1\rangle+\langle
u_2,u_2\rangle\\
&=&b^*b+\Phi_1(u_1^*u_1)+\Phi_2(u_2^*u_2)\\
&\geq&0
\end{eqnarray*}

  \end{proof}

\begin{theorem}\label{condpositivity}
 Let $\mathfrak{B}$ be a $C^*$-algebra and $\mathfrak{D}$  a $C^*$-subalgebra of $\mathfrak{B}$.
 Suppose that $\mathfrak{A}_1$, $\mathfrak{A}_2$  are $\ast$-algebras containing
 $\mathfrak{B}$,
 each endowed with two positive conditional
 expectations
 $\Phi_j:\mathfrak{A}_j\longrightarrow\mathfrak{B}$,  and $\Psi_j:\mathfrak{A}_j\longrightarrow\mathfrak{D}$,
 $j=1,2$ and consider the c-free product
 $(\mathfrak{A}, \Phi,\Psi)=\free_{i=1,2}(\mathfrak{A}_i,\Phi_i,\Psi_i)$.

Then the maps $\Phi$ and $\Psi$ are positive.
\end{theorem}\label{condposit}
\begin{proof}

The positivity of $\Psi$ is by now a classical result in the theory
of free probability with amalgamation over a $C^*$-algebra (for
example, see \cite{speicher1}, Theorem 3.5.6). For the positivity of
$\Phi$ we have to show that $\Phi(a^*a)\geq 0$ for any
$a\in\mathfrak{A}$.

Any element of $\mathfrak{A}$ can be written as
\begin{eqnarray*}
a&=& \sum_{k=1}^{N}s_{1,k}\dots s_{n(k),k},
\end{eqnarray*}
where $s_{j,k}\in\mathfrak{A}_{\varepsilon(j,k)}\
\varepsilon(1,k)\neq\varepsilon(2,k)\neq\dots\neq
\varepsilon(n(k),k).$

\noindent Writing
\begin{eqnarray*}
 s_{(j,k)}&=&
s_{(j,k)}-\Psi(s_{(j,k))}+\Psi(s_{(j,k))}
\end{eqnarray*}
and expanding the product, we can consider $a$ of the form
 \[
 a=d+\sum_{k=1}^{N}a_{1,k}\dots a_{n(k),k}
  \]
with
 $d\in\mathfrak{D}\subset \mathfrak{B}$
  and
$a_{j,k}\in\mathfrak{A}_{\varepsilon(j,k)}$
 such that
 $\Psi_{\varepsilon(j,k)}(a_{j,k})=0$
 and
$\varepsilon(1,k)\neq\varepsilon(2,k)\neq\dots\neq
\varepsilon(n(k),k)$.

Therefore
\begin{eqnarray*}
\Phi(a^*a)
 &=&
 \Phi\Big(d^*d+d^*
  \Big[
  \sum_{k=1}^{N}a_{1,k}\dots a_{n(k),k}
   \Big]
 +
 \Big[
 \sum_{k=1}^{N}a_{1,k}\dots a_{n(k),k}
 \Big]^*d+
 \\
&&
 \Big[\sum_{k=1}^{N}a_{1,k}\dots
a_{n(k),k}\Big]^*\Big[\sum_{k=1}^{N}a_{1,k}\dots
a_{n(k),k}\Big]\Big).\\
\end{eqnarray*}
Since $\Phi$ is a conditional expectation and
$d\in\mathfrak{D}\subset\mathfrak{B}$, the above equality becomes

\begin{eqnarray*}
 \Phi(a^*a)
  &=&
d^*d+\sum_{k=1}^{N}d^*\Phi( a_{1,k}\dots
a_{n(k),k})+\sum_{k=1}^{N}\Phi( a_{n(k),k}^*\dots a_{1,k}^*)d\\
&&+ \sum_{k,l=1}^{N}\Phi\big( a_{n(k),k}^*\dots a_{1,k}^*
a_{1,l}\dots
a_{n(l),l}\big).\\
\end{eqnarray*}
Using the definition of the conditionally free product with
amalgamation over $\mathfrak{B}$ and that
$\Psi_{\varepsilon(j,k)}(a_{j,k})=0$ for all $j,k$, one further has
\begin{eqnarray*}
\Phi(a^*a)
 &=&
 d^*d+\sum_{k=1}^{N}
 \Phi(d^* a_{1,k})\Phi(a_{2,k})\dots\Phi(a_{n(k),k})\\
 &&
  +\sum_{k=1}^{N}\Phi( a_{n(k),k})^*\dots \Phi(a_{2,k}^*)\Phi(a_{1,k}^*d)\\
 &&
  +\sum_{k,l=1}^{N}
   \left[
   \Phi(a_{n(k),k})^*\dots\Phi(a_{2,k}^*)
    \right]
  \Phi(a_{1,k}^* a_{1,l})\Phi(a_{2,l})\dots\Phi(a_{n(l),l})\\
 \end{eqnarray*}
that is
\begin{eqnarray*}
\Phi(a^*a)
 &=&
 d^*d+\sum_{k=1}^{N}\Phi(d^*
 a_{1,k})\left[\Phi(a_{2,k})\dots\Phi(a_{n(k),k})\right]\\
 &&
  +\sum_{k=1}^{N}\left[\Phi(a_{2,k})\dots \Phi( a_{n(k),k})\right]^*\Phi(a_{1,k}^*d)\\
 &&
  +\sum_{k,l=1}^{N}\left[\Phi(a_{2,k})\dots \Phi( a_{n(k),k})\right]^*
   \Phi(a_{1,k}^* a_{1,l})\left[\Phi(a_{2,l})\dots\Phi(a_{n(l),l})\right]\\
\end{eqnarray*}

 From Lemma \ref{mixprod},
 the matrix $S=\left(\Phi(a_{1,i}^*a_{1,j})_{i,j=1}^{N+1}\right)$
 is positive in $M_{N+1}(\mathfrak{B})$, therefore
 \[S=T^*T,\  \text{for some}\ T\in M_{N+1}(\mathfrak{B}).\]

 Denote now $a_{1,N+1}=d$ and
$v_k=\Phi(a_{2,k})\dots\Phi(a_{n(k),k})$.The identity for
$\Phi(a^*a)$ becomes:
 \begin{eqnarray*}
\Phi(a^*a)&=&(v_1,\dots,v_N,1)^*T^*T(v_1,\dots,v_N,1)\\
&\geq&0,\ \ \text{as claimed.}
 \end{eqnarray*}
\end{proof}

\begin{theorem}\label{assoc}
Assume that $\mathcal{I}=\bigcup_{j\in\mathfrak{J}}\mathcal{I}_j$ is
a partition of $\mathcal{I}$. Then:
\[\free_{j\in\mathfrak{J}}
\left(\free_{i\in\mathcal{I}_j}(\mathfrak{A}_i,\Phi_i,\Psi_i)\right)
=\free_{i\in\mathcal{I}}(\mathfrak{A}_i,\Phi_i,\Psi_i)\]
\end{theorem}
\begin{proof}
The proof is identical to the proofs of similar results in
\cite{mlotk1} and \cite{bls}. Consider
$a_i\in\mathfrak{A}_{\varepsilon(i)}, 1\leq i\leq m$ such that
$\varepsilon(1)\neq\varepsilon(2\neq\dots\neq\varepsilon(m)$ and
$\Psi_{\varepsilon(i)}(a_i)=0$. Let $1=i_0<i_1<\dots <i_k= m$ and
$\mathfrak{J}_l=\{\varepsilon(i), i_{l-1}\leq i<i_l\}$.

\noindent Since
\[\left(\free_{j\in\mathfrak{J}_l}\Psi_j\right)((a_{i_{l}-1}\cdots
a_{i_l}))=0,\]
it suffices to show that
\begin{equation*}
\Phi(a_1\cdots
a_m)=\prod_{l=1}^k\left[(\free_{(\Psi_j),j\in\mathfrak{J}_l}\Phi_j)(a_{i_{l}-1}\cdots
a_{i_l})\right].
\end{equation*}
But
\[
\Phi(a_1\cdots
a_m)=\Phi_{\varepsilon(1)}(a_1)\cdots\Phi_{\varepsilon(m)}(a_m)
\]
while, since $\Psi_{\varepsilon(i)}(a_i)=0$,
\[ (\free_{(\Psi_j),j\in\mathfrak{J}_l}\Phi_j)(a_{i_{l}-1}\cdots
a_{i_l})=\Phi_{{i_{l}-1}}(a_{i_{l}-1})\cdots\Phi_{{i_{l}}}(a_{i_{l}})\]
and the conclusion follows.
\end{proof}

\begin{definition}Let $\mathfrak{A}$ be an algebra (respectively a $\ast$-algebra),
$\mathfrak{B}$ a subalgebra ($\ast$-subalgebra) of $\mathfrak{A}$
and $\mathfrak{D}$ a subalgebra ($\ast$-subalgebra) of
$\mathfrak{B}$. Suppose $\mathfrak{A}$ is endowed with the
conditional expectations
$\Psi:\mathfrak{A}\longrightarrow\mathfrak{D}$ and
$\Phi:\mathfrak{A}\longrightarrow\mathfrak{D}$.
\begin{enumerate}
\item[(i)]The subalgebras ($\ast$-subalgebras)
$(\mathfrak{A}_i)_{i\in\mathcal{I}}$ of $\mathfrak{A}$ are said to
be \emph{c-free} with respect to $(\Phi, \Psi)$ if:
\begin{enumerate}
\item[(a)]$(\mathfrak{A}_i)_{i\in\mathcal{I}}$ are free with respect to
$\Psi$.
\item[(b)]if $a_i\in\mathfrak{A}_{\varepsilon(i)}, 1\leq i\leq m$,
are such that $\varepsilon(1)\neq\dots\neq\varepsilon(m)$ and
$\Psi(a_i)=0$,
 then
 $
 \Phi(a_1\cdots a_m)=\Phi(a_1)\cdots \Phi(a_m).
 $
\end{enumerate}
\item[(ii)] The elements $(X_i)_{i\in\mathcal{I}}$ of $\mathfrak{A}$ are
said to be c-free with respect to $(\Phi, \Psi)$ if the subalgebras
($\ast$-subalgebras) generated by $\mathfrak{B}$ and $X_i$ are
c-free with respect to $(\Phi, \Psi)$.
\end{enumerate}
\end{definition}

We will denote by $\mathfrak{B}\langle\xi\rangle$ the
non-commutative algebra of polynomials in the symbol $\xi$ and with
coefficients from $\mathfrak{B}$ (the coefficients do not commute
with the symbol $\xi$). If $\mathcal{I}$ is a family of indices,
$\mathfrak{B}\langle\{\xi_i\}_{i\in \mathcal{I}}\rangle$ will denote
the algebra of polynomials in the non-commuting variables
$\{\xi\}_{i\in \mathcal{I}}$ and with coefficients from
$\mathfrak{B}$. We will identify
 $\mathfrak{B}\langle\{\xi_i\}_{i\in \mathcal{I}}\rangle$
   with the free product with amalgamation over
$\mathfrak{B}$ of the family
$\{\mathfrak{B}\langle\xi_i\rangle\}_{i\in \mathcal{I}}$.

If $\mathfrak{A}$ is a $\ast$-algebra and $\mathfrak{B}$ is with the
$C^*$-algebra, $\mathfrak{B}\langle\xi\rangle$ will also be
considered with a $\ast$-algebra structure, by taking $\xi^*=\xi$.
If $X$ is a selfadjoint element from $\mathfrak{A}$, we define the
conditional expectations
  \[
  \Phi_X, \Psi_X:\mathfrak{B}\langle\xi\rangle\longrightarrow\mathfrak{B}
  \]
 given by
$ \Phi_X(f(\xi))=\Phi(f(X))$
  and
  $\Psi_X(f(\xi))=\Psi(f(X))$ ,
for
any $f(\xi)\in\mathfrak{B}\langle\xi\rangle$.

\begin{corollary}\label{sumposit}
Suppose that $\mathfrak{A}$ is a $\ast$-algebra and $X$ and $Y$ are
c-free selfadjoint elements of $\mathfrak{A}$ such that the maps
$\Phi_X, \Psi_X$ and $\Phi_Y, \Psi_Y$ are positive. Then the maps
$\Phi_{X+Y}$ and $\Psi_{X+Y}$ are also positive.
\end{corollary}
\begin{proof}
The positivity of $\Psi_{X+Y}$ is an immediate consequence of the
fact that $X$ and $Y$ are free with amalgamation over $\mathfrak{B}$
with respect to $\Psi$. It remains to prove the positivity of
$\Phi_{X+Y}$.

 Since the maps $\Phi_X:\mathfrak{B}\langle\xi_1\rangle\longrightarrow\mathfrak{B}$ and
$\Phi_Y:\mathfrak{B}\langle\xi_2\rangle\longrightarrow\mathfrak{B}$
are positive, from Theorem \ref{condpositivity} so is
 \[
\Phi_x\ast_{(\Psi_X,\Psi_Y)}\Phi_Y:\mathfrak{B}\langle\xi_1\rangle\ast_\mathfrak{B}
\mathfrak{B}\langle\xi_2\rangle
=\mathfrak{B}\langle\xi_1,\xi_2\rangle\longrightarrow\mathfrak{B}
 \]

Remark also that
\[i_Z:\mathfrak{B}\langle\xi\rangle\ni f(\xi)\mapsto f(X+Y)
\in\mathfrak{B}\langle\xi_1\rangle\ast_\mathfrak{B}\mathfrak{B}\langle\xi_2\rangle\]
is a positive $\mathfrak{B}$-functional.

The conclusion follows from the fact that the c-freeness of $X$ and
$Y$ is equivalent to
\[\Phi_{X+Y}=(\Phi_X\ast_{(\Psi_X,\Psi_Y)}\Phi_Y)\circ i_{X+Y}.\]

\end{proof}

%%%%%%%%%%%%%%%%%%%%%%%%%%%%%%%%%%%%%%%%%%%%%%%%%%%%%%%%%%%%%%%%%%%%%%%%%%%%%%%%%%%%%%%%%%%%%%%%%%%%%%%
%%%%%%%%%%%333333333333333333333333333333333333333333333333333333333333333333333333333333

\section{Multilinear function series and the $^cR$-transform}

 Let $\mathfrak{A}$ be a $\ast$-algebra containing the $C^*$-algebra $\mathfrak{B}$,
 endowed with a conditional expectation $\Psi:\mathfrak{A}\longrightarrow\mathfrak{B}$.  If $X$ is
 a selfadjoint element of $\mathfrak{A}$, then by the
  moment of order $n$ of $X$ we will understand the map
\begin{eqnarray*}
m_X^{(n)}
 &:&
  \underbrace{\mathfrak{B}\times\dots\times\mathfrak{B}}_{n-1\ \text{times}}\longrightarrow\mathfrak{B} \\
m_X^{(n)}(b_1,\dots,b_{n-1})
 &=&
 \Psi(Xb_1X\dots Xb_{n-1}X)
\end{eqnarray*}

If $\mathfrak{B}=\mathbb{C}$, then the moment-generating series of
$X$
\[m_X(z)=\sum_{n=0}^\infty\Psi(X^n)z^n\]
encodes all the information about the moments of $X$. For
$\mathfrak{B}\neq \mathbb{C}$, the straightforward generalization
\[\overline{m}_X(z)=\sum_{n=0}^\infty\Psi(X^n)z^n\]
generally fails to keep track of all the possible moments of $X$. A
solution to this inconvenience was proposed in \cite{dykema}, namely
the moment-generating multilinear function series of $X$. Before
defining this notion, we will briefly recall the construction and
several results on multilinear function series.

Let $\mathfrak{B}$ be an algebra over a field $K$. We set
$\widetilde{\mathfrak{B}}$ equal to $\mathfrak{B}$ if $\mathfrak{B}$
is unital and to the unitization of $\mathfrak{B}$ otherwise. For
$n\geq1$, we denote by $\mathcal{L}_n(\mathfrak{B})$ the set of all
$K$-multilinear maps
 \[
\omega_n:\underbrace{\mathfrak{B}
\times\dots\times\mathfrak{B}}_{n\ \text{times}}
\longrightarrow\mathfrak{B}
 \]

A \emph{formal multilinear function series} over $\mathfrak{B}$ is a
sequence $\omega=(\omega_0,\omega_1,\dots)$, where
$\omega_0\in\widetilde{\mathfrak{B}}$ and
$\omega_n\in\mathcal{L}_n(\mathfrak{B})$ for $n\geq1$. According to
\cite{dykema}, the set of all multilinear function series over
$\mathfrak{B}$ will de denoted by $\mul$.

For $\alpha,\beta\in\mul$, the \emph{ formal sum} $\alpha+\beta$ and
the \emph{formal product} $\alpha\beta$ are the elements from $\mul$
defined by:

\begin{eqnarray*}
(\alpha+\beta)_n(b_1,\dots,b_n)&=&\alpha_n(b_1,\dots,b_n)+\beta_n(b_1,\dots,b_n)\\
(\alpha\beta)_n(b_1,\dots,b_n)&=&
\sum_{k=0}^n\alpha_k(b_1,\dots,b_k)\beta_{n-k}(b_{k+1},\dots,b_n)\\
\end{eqnarray*}
for any $b_1,\dots,b_n\in\mathfrak{B}$.

If $\beta_0=0$, then the \emph{formal composition} $\alpha\circ
\beta\in\mul$ is defined by
\[
 (\alpha\circ\beta)_0=\alpha_0 \hspace{4.7cm}
  \]
and, for $n\geq1$, by
\begin{eqnarray*}
(\alpha\circ \beta)_n(b_1,\dots,b_n)&=&\sum_{k=1}^n
%\hspace{-.2cm}
\sum_{
%\substack{p_1,\dots,p_k\geq1\\
%{p_1+\dots+p_k=n} }
}
%\hspace{-.5cm}
\alpha_k\Big(\beta_{p_1}(b_1,\dots,b_{p_1}),\dots,
\beta_{p_k}(b_{q_k+1},\dots,b_{q_k+p_k})\Big)
\end{eqnarray*}
\noindent where  the second summation is done over all $k$-tuples
$p_1,\dots,p_k\geq1$ such that $p_1+\dots+p_k=n$ and
$q_j=p_1+\dots+p_{j-1}$.

One can work with elements of $\mul$ as if they were formal power
series. The relevant properties are described in  \cite{dykema},
Proposition 2.3 and Proposition 2.6. As in \cite{dykema}, we use
$1$, respectively $I$, to denote the identity elements of $\mul$
relative to multiplication, respectively composition. In other
words, $1=(1,0,0,\dots)$ and $I=(0, id_\mathfrak{B},0,0,\dots)$. We
will also use the fact that an element $\alpha\in\mul$ has an
inverse with respect to formal composition, denoted
$\alpha^{{\langle -1\rangle}}$, if and only if $\alpha$ has the form
$(0,\alpha_1,\alpha_2,\dots)$ with $\alpha_1$ an invertible element
of $\mathcal{L}_1(\mathfrak{B})$.

\begin{definition}
\emph{With the above notation, the moment-generating multilinear
function series $\mathcal{M}_X$ of $X$ is the element of $\mul$ such
that:}
\begin{eqnarray*}
\mathcal{M}_{X,0}&=&\Psi(X)\\
\mathcal{M}_{X,n}(b_1,\dots,b_n)&=&\Psi(Xb_1X\cdots Xb_nX).
\end{eqnarray*}
\end{definition}

Given an element $\alpha\in\mul$, the multilinear function series
$R_\alpha$ is defined by the following equation (see \cite{dykema},
Def 6.1):

\begin{equation}\label{ralpha}
R_\alpha=\left(\left(1+\alpha I\right)^{-1}\right)\circ(I+I\alpha
I)^{\langle-1\rangle}.
\end{equation}

A key property of $R$ is that for any $X,Y\in\mathfrak{A}$ free over
$\mathfrak{B}$, we have

\begin{equation}\label{ralpha1}
R_{\mathcal{M}_{X+Y}}=R_{\mathcal{M}_{X}}+R_{\mathcal{M}_{Y}}.
\end{equation}

These relations were proved earlier in the particular case
$\mathfrak{B}=\mathbb{C}$. One can also describe $R_\alpha$ by
combinatorial means, via the recurrence relation

\begin{eqnarray*}
\alpha_n(b_1,\dots,b_n)
 &=&
 \sum_{k=0}^{n} \sum
R_{\alpha,k}\Big([b_{1}\alpha_{p(1)}(b_3,\dots,b_{i_1-2})b_{i_1-1}], \dots\\
&&
 \hspace{-1.6cm}\dots, [b_{i(k-1)}\alpha_{p(k)}(b_{i(k-1)+1},
\dots,b_{i(k)-2})b_{i(k)-1}]\Big)b_{i(k)}
\alpha_{n-i_k}(b_{i_{k+1}},\dots,b_n)
\end{eqnarray*}
where the second summation is done over all
$1=i(0)<i(1)<\dots<i(k)\leq n$ and $p(k)=i(k)-i(k-1)-2$.

Following an idea from \cite{bls}, the above equation can be
graphically illustrated by the picture:
 \setlength{\unitlength}{.16cm}
 \begin{equation*}
 \begin{picture}(10,8)

 \put(-25,2){\line(0,1){3}}
 \put(-25,2){\line(1,0){10}}
\put(-15,2){\line(0,1){3}} \put(-25,5){\line(1,0){10}}
\put(-14,3){=}

\put(-12,3){\Huge{$\sum$}}

\put(-7,1){\line(0,1){4}}

\put(-7,1){\line(1,0){14}}

\put(-1,1){\line(0,1){4}}

\put(5,1){\line(0,1){4}}

\put(21,1){\line(0,1){4}}

\put(21,1){\line(-1,0){7}}

\put(15,1){\line(0,1){4}} \qbezier[7](7,1)(9,1)(13,1)

%moments
%\put(-6,2){\rule{4\unitlength}{3\unitlength}}

\put(-6,2){\line(0,1){3}} \put(-2,2){\line(0,1){3}}

\put(-6,2){\line(1,0){4}}

\put(-6,5){\line(1,0){4}}

\put(0,2){\line(0,1){3}} \put(4,2){\line(0,1){3}}

\put(0,2){\line(1,0){4}}

\put(0,5){\line(1,0){4}}

\put(16,2){\line(0,1){3}} \put(20,2){\line(0,1){3}}

\put(16,2){\line(1,0){4}}

\put(16,5){\line(1,0){4}}

\put(23,2){\line(0,1){3}}

\put(23,2){\line(1,0){9}}

\put(32,2){\line(0,1){3}}

\put(23,5){\line(1,0){9}}

 \end{picture}
 \end{equation*}

In the case of scalar c-free probability, an analogue of the
Voiculescu's $R$-transform is developed in \cite{bls}. In order to
avoid confusions, we will denote it by ${}^cR$.

The ${}^cR$-transform has the property that it linearizes the c-free
convolution of pairs of compactly supported measures. In particular,
if $X$ and $Y$ are c-free elements from some algebra $\mathfrak{A}$,
then
\[{}^cR_{X+Y}={}^cR_{X}+{}^cR_{Y}.\]

If the $\ast$-algebra $\mathfrak{A}$ is endowed with the
$\mathbb{C}$-valued states $\varphi, \psi$ and $X$ is a selfadjoint
element of $\mathcal{A}$, then (see \cite{bls}), the coefficients
$\{{}^cR_m\}_m\geq0$ of ${}^cR_X$ are defined by the recurrence:
\begin{eqnarray*}
\varphi(X^n)
 &=&
  \sum_{k=1}^n
  \sum_ {
  \substack{
l(1),\dots,l(k)\geq 0\\
 l(1)+\dots+l(k)=n-k
 }}
  {}^cR_k\cdot
 \psi(X^{l(1)})\cdots\psi(X^{l(k-1)})\varphi(X^{l(k)})
\end{eqnarray*}
equation that can be graphically illustrated by the picture, were
the dark boxes stand for the application of $\varphi$ and the light
ones for the application of $\psi$:

 \setlength{\unitlength}{.16cm}
 \begin{equation*}
 \begin{picture}(10,6)

 \put(-25,2){\rule{10\unitlength}{3\unitlength}}
 \put(-14,3){=}

\put(-12,3){\Huge{$\sum$}}

\put(-7,1){\line(0,1){4}}

\put(-7,1){\line(1,0){14}}

\put(-1,1){\line(0,1){4}}

\put(5,1){\line(0,1){4}}

\put(21,1){\line(0,1){4}}

\put(21,1){\line(-1,0){7}}

\put(15,1){\line(0,1){4}} \qbezier[7](7,1)(9,1)(13,1)

%moments
%\put(-6,2){\rule{4\unitlength}{3\unitlength}}

\put(-6,2){\line(0,1){3}} \put(-2,2){\line(0,1){3}}

\put(-6,2){\line(1,0){4}}

\put(-6,5){\line(1,0){4}}

\put(0,2){\line(0,1){3}} \put(4,2){\line(0,1){3}}

\put(0,2){\line(1,0){4}}

\put(0,5){\line(1,0){4}}

\put(16,2){\line(0,1){3}} \put(20,2){\line(0,1){3}}

\put(16,2){\line(1,0){4}}

\put(16,5){\line(1,0){4}}

\put(23,2){\rule{9\unitlength}{3\unitlength}}

 \end{picture}
 \end{equation*}
The above considerations lead to the following definition:
\begin{definition}\label{crcar}
Let $\beta, \gamma\in\mul$. The multilinear function series
${}^cR_{\beta,\gamma}$ is the element of $\mul$ defined by the
recurrence relation
\begin{eqnarray*}
\beta_n(b_1,\dots,b_n)
 &=&
   \sum_{k=0}^{n} \sum
{}^cR_{\beta, \gamma,k}
 \Big([b_{1}\gamma_{p(1)}(b_3,\dots,b_{i_1-2})b_{i_1-1}], \dots\\
&&
 \hspace{-1.4cm}\dots, [b_{i(k-1)}\gamma_{p(k)}(b_{i(k-1)+1},
\dots,b_{i(k)-2})b_{i(k)-1}]
 \Big)b_{i(k)}
\beta_{n-i_k}(b_{i_{k+1}},\dots,b_n)
\end{eqnarray*}
where the second summation is done over all
$1=i(0)<i(1)<\dots<i(k)\leq n$ and $p(k)=i(k)-i(k-1)-2$.
\end{definition}

The following analytical description of ${}^cR_{\beta, \gamma}$ also
shows that it is unique and well-defined:
\begin{theorem}\label{cranalit} For any $\beta, \gamma\in\mul$,
\begin{equation}\label{analit2}
^cR_{\beta,\gamma}=\left[\beta(1+I\beta)^{-1}\right]\circ(I+I\gamma
I)^{\langle-1\rangle}
\end{equation}
\end{theorem}
Before proving the theorem, remark that the right-hand side  of
(\ref{analit2}) is well-defined and unique, since $1+I\gamma$ is
invertible with respect to the formal multiplication, $I+I\beta I$
is invertible with respect to formal composition and its inverse has
0 as first component (see \cite{dykema}). We will need the following
auxiliary result:

 \begin{lemma}\label{aux1}
 Let $\beta$ be an element of $\mul$ and $I$ the identity element
 with respect to formal composition, $I=(0,id_\mathfrak{B}, 0,0\dots)$.
 \begin{enumerate}
 \item[(i)]the multilinear function series $I\beta$ is given by:
 \begin{eqnarray*}
 (I\beta)_0&=&0\\
 (I\beta)_n(b_1,\dots,b_n)&=&b_1\beta_{n-1}(b_2,\dots,b_n)
 \end{eqnarray*}
 \item[(ii)]the multilinear function series $I\beta I$ is given by
  \begin{eqnarray*}
(I\beta I)_0&=&0\\
(I\beta I)_1(b_1)&=&0\\
 (I\beta I)_n(b_1,\dots,b_n)&=&b_1\beta_{n-2}(b_2,\dots,b_{n-1})b_n
\end{eqnarray*}
 \end{enumerate}
 \end{lemma}

\begin{proof}
Since $I=(0, id_\mathfrak{B}, 0,\dots)$, one has:

\[(I\beta)_0=I_0\beta_0=0.\]
If $n\geq1$,
\begin{eqnarray*}
(I\beta)_n(b_1, \dots, b_n)&=&\sum_{k=0}^n I_k(b_1, \dots,
b_k)\beta_{n-k}(b_{k+1},\dots, b_n)\\
&=&I_1(b_1)\beta_{n-1}(b_{k+1},\dots, b_n)\\
&=&b_1\beta_{n-1}(b_{k+1},\dots, b_n).
\end{eqnarray*}

For $I\beta I$, the same computations give:

\begin{eqnarray*}
(I\beta I)_0&=&(I \beta)_0 I_0=0\\
(I\beta I)_1&=&(I \beta)_0 I_1(b_1)+(I \beta)_1(b_1) I_0\\
&=&0.
\end{eqnarray*}
If $n\geq2$, one has:
\begin{eqnarray*}
(I\beta I)_n(b_1,\dots,b_n) &=&\sum_{k=0}^n(I\beta)_k(b_1,\dots,
b_k)I_{n-k}(b_{k+1},\dots,b_n)\\
&=&(I\beta)_{n-1}(b_1,\dots,b_k)I_1(b_1)\\
 &=&b_1\beta_{n-2}(b_2,\dots,b_{n-1})b_n
\end{eqnarray*}
\end{proof}

\begin{proof}[Proof of the Theorem \ref{cranalit}:]

Set $\sigma=I+I\beta I$. Then
\begin{eqnarray*}
\left({}^cR_{\beta,\gamma}\circ\sigma\right)_n(b_1,\dots,b_n)&=&\sum_{k=1}^n
%\hspace{-.2cm}
\sum_{
\substack{p_1,\dots,p_k\geq1\\
{p_1+\dots+p_k=n}
 } } \hspace{-.5cm}
{}^cR_{\beta,\gamma,k}\Big(\sigma_{p_1}(b_1,\dots,b_{p_1}),\dots,\\
&&
 \hspace{3.7cm}\sigma_{p_k}(b_{q_k+1},\dots,b_{q_k+p_k})\Big)
\end{eqnarray*}
where $q_i=p_1+\dots+p_{i-1}$.

From Lemma (\ref{aux1})(ii), we have that

\[
 \sigma_n(b_1,\dots,b_n)=(I+I\beta I)_n(b_1,\dots,b_n)
 \]
therefore Definition \ref{crcar} is equivalent to
\[
 \beta_n(b_1,\dots,b_n)
  =
\sum_{k=0}^{n}
 \left(
   {}^cR_{\beta,\gamma}\circ(I+I\beta
I)_{k}(b_1,\dots,b_k)
  \right)
   b_{k+1} \beta_{n-k-2}(b_{k+2},\dots,b_n)
\]
Considering now Lemma \ref{aux1}(i), the above relation becomes
\[
 \beta_n(b_1,\dots,b_n)
  =
\sum_{k=0}^{n}
 \left(
  {}^cR_{\beta,\gamma}\circ(I+I\beta
I)_{k}(b_1,\dots,b_k)
 \right)
  (I+I\beta )_{n-k}(b_{k+1} ,\dots,b_n)
\]
therefore
\[
 \beta=\left[^cR_{\beta,\gamma}\circ(I+I\gamma I)\right](1+I\beta)
 \]

which is equivalent to (\ref{analit2}).
\end{proof}

\begin{remark}
Up to a shift in the coefficients,  equation (\ref{analit2}) is
similar to the result in the case $\mathfrak{B}=\mathbb{C}$ from
\cite{bls}, Theorem 5.1.
\end{remark}

Let $X$ be a selfadjoint element of $\mathfrak{A}$. If
$\mathfrak{A}$ is endowed with two $\mathfrak{B}$-valued conditional
expectations $\Phi, \Psi$, the element $X$ will have two
moment-generating multilinear function series, one with respect to
$\Psi$, that we will denote by $\mathcal{M}_X$, and one with respect
to $\Phi$, denoted $\mathfrak{M}_X$. For brevity, we will use the
notation $^cR_X$ for the multilinear function series
$^cR_{\mathcal{M}_X, \mathfrak{M}_X}$.
\begin{theorem}\label{adit}
Let $X$ and $Y$ be two elements of $\mathfrak{A}$ that are c-free
with respect to the pair of conditional expectations $(\Phi, \Psi)$.
Then
\[^cR_{X+Y}={}^cR_X+{}^cR_Y\]
\end{theorem}
\begin{proof}
 Let $\mathcal{A}$ be an algebra containing $\mathfrak{B}$ as a subalgebra
 and endowed with the conditional expectations $\Phi,
 \Psi:\mathcal{A}\longrightarrow\mathfrak{B}$. Consider the set
 $\mathcal{A}_0=\mathcal{A}\setminus\mathfrak{B}$ (set difference). For
 $n\geq1$ define the maps
 \[
  {}^cr:\underbrace{\mathcal{A}_0\times\dots\times\mathcal{A}_0}_{n\ \text{times}}
 \longrightarrow\mathfrak{B}
  \]
 given by the recurrence formula:
 \begin{eqnarray*}
 \Phi(a_1\cdots a_n)&=&\sum_{k=1}^n
 \sum_{
\substack{l(1)<\dots<l(k)\\
1<l(1), l(k)\leq n\\
}}
 {}^cr_{k}\text{\huge{$($}}a_1[\Psi(a_2\cdots a_{l(1)-1})],\dots,\\
  &&\dots,a_{l(k-1)}[\Psi(a_{l(k-1)+1}\cdots
  a_{l(k)_1})],a_{l(k)}[\Phi(a_{l(k)+1}\cdots a_n)]
 \text{\huge{$)$}}
 \end{eqnarray*}

 Note that ${}^cr_n$ is well defined, and that, for any
 $b_1,\dots,b_n\in\mathfrak{B}$,
 \begin{equation}\label{simi1}
 {}^cr_{n+1}(X,b_1X,\dots ,b_nX)={}^cR_{X,n}(b_1,\dots,b_n).
 \end{equation}

 As in Section \ref{second},  consider $\mathfrak{B}\langle\xi_i\rangle$,
 the noncommutative algebras of polynomials in
 the symbols $\xi_i, i=1,2$ and with coefficients from $\mathfrak{B}$ and the
  conditional
 expectations
  \[
   \Phi_X, \Psi_X:\mathfrak{B}\langle\xi_1\rangle\longrightarrow\mathfrak{B}
    \]
 given by
 \begin{eqnarray*}
\Phi_X(f(\xi_1))&=&\Phi(f(X))\\
\Psi_X(f(\xi_1))&=&\Psi(f(X))
 \end{eqnarray*}
and their analogues $\Phi_Y, \Psi_Y$ for
$\mathfrak{B}\langle\xi_2\rangle$.

On $\mathfrak{B}\langle\xi_1,\xi_2\rangle$, identified to
$\mathfrak{B}\langle\xi_1\rangle\ast_{\mathfrak{B}}\mathfrak{B}\langle\xi_2\rangle$,
consider the conditional expectations $\Psi_0, \Phi_0,\varphi$ given
by:
\begin{eqnarray*}
\Psi_0&=&
\Psi_X\ast\Psi_Y\\
\Phi_0(f(\xi_1,\xi_2))
 &=&
  \Phi(f(X,Y))\\
\varphi(a_1a_2\dots a_n)
 &=&
   \sum_{k=1}^n\sum_
{\substack{l(1)<\dots<l(k)\\
1<l(1), l(k)\leq n\\
}}
\rho_k\text{\huge{$($}}a_1[\Psi_0(a_2\cdots a_{l(1)-1})],\dots,\\
&&\dots,
  a_{l(k-1)}[\Psi_0(a_{l(k-1)+1}\cdots
  a_{l(k)_1})],a_{l(k)}[\varphi(a_{l(k)+1}\cdots a_n)]
  \text{\huge{$)$}}\\
\end{eqnarray*}
where $a_1,\dots,a_n$ are elements of the set
$
 \mathfrak{B}\langle\xi_1,\xi_2\rangle_0
=\mathfrak{B}\langle\xi_1\rangle\cup\mathfrak{B}\langle\xi_2\rangle\setminus\mathfrak{B},
 $
 and  the maps
\[\rho_n:\underbrace{\mathfrak{B}\langle\xi_1,\xi_2\rangle_0\times
\dots\mathfrak{B}\langle\xi_1,\xi_2\rangle_0} _{n\
\text{times}}\longrightarrow\mathfrak{B}\]
are given by:
\[
 \rho_n(a_1,\dots,a_n)=\left\{\begin{array}{cl}
{}^cr(a_1,\dots,a_n)& \text{if all}\ a_1,\dots, a_n\in\mathfrak{B}\langle\xi_1\rangle\\
{}^cr(a_1,\dots,a_n)& \text{if all}\ a_1,\dots, a_n\in\mathfrak{B}\langle\xi_2\rangle\\
0& \text{otherwise} \\
\end{array}\right..
 \]

We will show that $\varphi=\Phi_0$, in particular $\varphi$ is also
well-defined. Consider the element
$a\in\mathfrak{B}\langle\xi_1,\xi_2\rangle$ of the form $a=a_1\cdots
a_n$ with $a_j\in\mathfrak{B}\langle\xi_{\varepsilon(j)}\rangle,$
such that
$\varepsilon(1)\neq\varepsilon(2)\neq\dots\neq\varepsilon(n)$ and
$\Psi_0(a_j)=0$. The computation of $\varphi(a_1\cdots a_n)$ is done
via the recurrence relation above. Because of the definition of
$\rho$ and the fact that $\Psi_0=\Psi_X\ast\Psi_Y$, only the term
with $k=1$ contribute at the sum, i.e.
\begin{eqnarray*}
\varphi(a_1\cdots a_n)&=&\varphi(a_1\varphi(a_2\cdots a_n))\\
&=&\varphi_{\varepsilon(1)}(a_1\varphi(a_2\cdots a_n)\\
&=&\varphi_{\varepsilon(1)}(a_1)\varphi(a_2\cdots a_n)
\end{eqnarray*}
and the identity between $\varphi$ and $\Phi_0$ follows by induction
over $n$.

Since $\varphi=\Phi_0$, the maps $\rho_n$ and ${}^cr_n$ are
satisfying the same recurrence relation, hence
\[\rho_n(a_1,\dots,a_n)={}^cr(a_1,\dots,a_n).\]
In particular
\begin{eqnarray*}
{}^cR_{X+Y,n}(b_1,\dots,b_n) &=&
{}^cr_{n+1}((X+Y)b_1(X+Y)\dots(X+Y)b_n(X+Y))\\
&&\hspace{-1.7cm}=
\rho_{n+1}((X+Y)b_1(X+Y)\dots(X+Y)b_n(X+Y))\\
&&\hspace{-1.7cm}=
\rho_{n+1}((X)b_1(X)\dots(X)b_n(X))+\rho_{n+1}((Y)b_1(Y)\dots(Y)b_n(Y))\\
&&\hspace{-1.7cm}=
{}^cR_{X,n}(b_1,\dots,b_n)+{}^cR_{Y,n}(b_1,\dots,b_n).
\end{eqnarray*}

\end{proof}

%%%%%%%%%%%%%%%%%%%%%%%%%%%%%%%%%%%%%%%%%%%%%%%%%%%%%%%%%%%%%%%%%%%%%%%%%%%%%%%%%%%%%%
%%%%%%%%%%%%%%%%%%%%%%%%44444444444444444444444444444444444444444444444

\section{Central limit theorem}

Consider the ordered set $\langle n \rangle=\{1,2,\dots, n\}$ and
$\pi$ a partition of $\langle n \rangle$ with blocks $B_1,\dots,
B_m$:
\[
\langle n\rangle=B_1\sqcup B_2\sqcup\dots\sqcup B_m.
\]

The blocks $B_p$ and $B_q$ of $\pi$ are said to be \emph{crossing}
if there exist $i<j<k<l$ in $\langle n\rangle$ such that $i,k\in
B_p$ and $j,l\in B_q$.

The partition $\pi$ is said to be \emph{non-crossing} if all pairs
of distinct blocks of $\pi$ are not crossing. We will denote by
$NC_2(n)$ the set of all non-crossing partitions of $\langle
n\rangle$ whose blocks contain exactly 2 elements and by $NC_{\leq
s}(n)$ the set of all non-crossing partitions of $\langle n\rangle$
whose blocks contain at most $s$ elements.

Let now $\gamma$ be a non-crossing partition of $\langle n\rangle$
and $B$ and $C$ be two blocks of $\pi$. We say that $B$ is interior
to $C$ if there exist two indices $i<j$ in $\langle n\rangle$ such
that $i,j\in C$ and $B\subset\{i+1,\dots,j-1\}$. The block $B$ is
said to be \emph{outer} if it is not interior to any other block of
$\gamma$. In a non-crossing partition of  $\langle n\rangle$, the
block containing $1$ is always outer.

Consider now an element $X$ of $\mathfrak{A}$. Let $\pi$ be a
partition from $NC_2(n+1)$ ($n$ = odd) and $B_1=(1,k)$ be the block
of $\pi$ containing $1$. We define, by recurrence, the following
expressions:

\begin{eqnarray*}
V_\pi(X,b_1,\dots,b_n) &=&
 \Psi(Xb_1 V_{\pi|\{2,\dots,j-1\}}
(X,b_2,\dots,b_{k-2})b_{k-1}X)b_k\\
&&
 \hspace{2cm}V_{\pi|\{k+1,\dots,n+1\}}(X,b_{k+1},\dots,b_n)\\
W_\pi(X,b_1,\dots,b_n)
 &=&
  \Phi(Xb_1 V_{\pi|\{2,\dots,j-1\}}
(X,b_2,\dots,b_{k-2})b_{k-1}X)b_k\\
&&
 \hspace{2cm}W_{\pi|\{k+1,\dots,n+1\}}(X,b_{k+1},\dots,b_n)\\
\end{eqnarray*}

\begin{theorem}\label{clthm}(Central Limit Theorem)
Let $(X_n)_{n\geq1}$ be a sequence of c-free elements of
$\mathfrak{A}$ such that:
\begin{enumerate}
\item[(1)]all $X_n$ have the same moment-generating multilinear
function series, $\mathfrak{M}$ with respect to $\Phi$ and $M$ with
respect to $\Psi$.
\item[(2)]$\Psi(X_n)=\Phi(X_n)=0.$
\end{enumerate}
Set
\[S_N=\frac{X_1+\dots+X_N}{\sqrt{N}},\]
Then:
\begin{enumerate}
\item[(i)]$\displaystyle
\lim_{N\rightarrow\infty}{}^cR_{S_N}=
(0,\mathfrak{M}_1(\cdot),0,\dots)$
\item[(ii)]$\displaystyle
\lim_{N\rightarrow\infty}R_{S_N}= (0,M_1(\cdot),0,\dots)$
\item[(iii)]there exist two conditional expectations
$\nu:\mathfrak{B}\langle\xi\rangle\longrightarrow\mathfrak{B}$,
depending only on $M_1(\cdot)$, and
$\mu:\mathfrak{B}\langle\xi\rangle\longrightarrow\mathfrak{B}$,
depending only on $M_1(\cdot)$ and $\mathfrak{M}_1(\cdot)$, such
that
\begin{eqnarray*}
\lim_{N\rightarrow\infty}\Psi_{S_N}&=&\nu\\
\lim_{N\rightarrow\infty}\Phi_{S_N}&=&\mu
\end{eqnarray*}
in the weak sense; in particular,
\begin{eqnarray*}
\nu(\xi b_1\xi\dots b_n\xi)&=&\sum_{\pi\in NC_2(n)}V_\pi(X_1,b_1,\dots,b_n)\\
\mu(\xi b_1\xi\dots b_n\xi)&=&\sum_{\pi\in
NC_2(n)}W_\pi(X_1,b_1,\dots,b_n).
\end{eqnarray*}
\end{enumerate}
\end{theorem}

\begin{proof}
Let $X$ be an element of $\mathfrak{A}$ with the same moment
generating series as $X_j,$\  $j\geq1$. As shown in \cite{dykema},
\[R_{S_N}=\sum_{k=1}^NR_{\frac{X_k}{\sqrt{N}}}=N R_{\frac{X}{\sqrt{N}}}. \]
Also, from Theorem \ref{assoc} and Theorem \ref{adit}, it follows
that
\[{}^cR_{S_N}=\sum_{k=1}^N{}^cR_{\frac{X_k}{\sqrt{N}}}=N {}^cR_{\frac{X}{\sqrt{N}}}. \]

Since $R$ and ${}^cR$ are multilinear and $M_0=\mathfrak{M}_0=0$, we
have that
\begin{eqnarray*}
\lim_{N\rightarrow\infty}{}^cR_{S_N,n}&=&
\lim_{N\rightarrow\infty}\frac{N}{N^{\frac{n+1}{2}}}{}^cR_{X,n}\\
&=&\left\{
\begin{matrix}
0 & \text{if}\ n\neq 1\\
\hspace{.3cm}\\
\mathfrak{M}_1(\cdot) & \text{if}\ n=1\\
\end{matrix}\right.
\end{eqnarray*}
and the similar relations for $R_{S_N,n}$, hence (i) and (ii) are
proved.

For (iii) it suffices to check the relations for $\nu(\xi
b_1\xi\dots b_n\xi)$ and $\mu(\xi b_1\xi\dots b_n\xi)$, which are a
trivial corollary of (i), (ii), and the recurrence formulas that
define $R$ and ${}^cR$.
\end{proof}
\begin{remark}\emph{For $\mathfrak{B}=\mathbb{C}$, the theorem is a weaker version of
Theorem 4.3 from \cite{bls}. If $\Psi$ is $\mathbb{C}$-valued, then
the result is similar to Corollary 5.1  from \cite{mlotk1}. Also,
under the assumptions that for some $a,b\in\mathfrak{B}$ we have
that:
 \begin{eqnarray*}
 \lim_{N\rightarrow\infty}N\Psi(X_1\cdots X_N)&=&a\\
 \lim_{N\rightarrow\infty}N\Psi(X_1\cdots X_N)&=&b
 \end{eqnarray*}
the same techniques lead to a  Poisson-type limit Theorem, similar
to Corollary 2, Section 5 of \cite{mlotk1}.}
\end{remark}

In the following remaining pages we will describe the positivity of
the limit functionals $\mu$ and $\nu$ in terms of $\Phi$ and $\Psi$.
The central result is Corollary \ref{corposit}.

For simplicity, suppose that $\mathfrak{B}$ is a unital
$\ast$-algebra (otherwise, we can replace $\mathfrak{B}$ by its
unitisation). Consider the symbol $\xi$, the $\ast$-algebra
$\mathfrak{B}\langle\xi\rangle$ of polynomials in $\xi$ with
coefficients from $\mathfrak{B}$, as defined before, and consider
also the linear space $\mathfrak{B}\xi\mathfrak{B}$ generated by the
set $\{b_1\xi b_2;\ b_1,b_2\in\mathfrak{B}\}$ with the
$\mathfrak{B}$-bimodule structure given by
\[a_1b_1\xi b_2a_2=(a_1b_1)\xi(b_2a_2)\]
for all $a_1,a_2,b_1,b_2\in\mathfrak{B}$.

\begin{lemma}\label{ext1}
For any positive $\mathfrak{B}$-sesquilinear  pairing
$\langle\cdot,\cdot\rangle$
 on $\mathfrak{B}\xi\mathfrak{B}$ there exists a positive
conditional expectation $
 \varphi:\mathfrak{B}\langle\xi\rangle\longrightarrow\mathfrak{B}
  $
such that for any $b_1,b_2\in\mathfrak{B}$ one has that
\[
 \varphi(\xi b_1^*b_2\xi)=\langle b_1\xi,b_2\xi\rangle
  \]
\end{lemma}
\begin{proof}
Without loss of generality, we can suppose that $\mathfrak{B}$ is
unital (otherwise we can replace $\mathfrak{B}$ by its unitization).

 Consider the
Full Fock bimodule over $\mathfrak{B}\xi\mathfrak{B}$
\[
 \mathcal{F}\langle\xi\rangle=\mathfrak{B}\oplus\Big(
\bigoplus_{n\geq
1}\underbrace{\mathfrak{B}\xi\mathfrak{B}\otimes_\mathfrak{B}\dots\otimes_\mathfrak{B}\mathfrak{B}\xi\mathfrak{B}}
_{n\ \text{times}}\Big)
 \]
with the pairing given by
\begin{eqnarray*}
\langle a,b\rangle&=& a^*b\\
\langle a_1\xi\otimes \dots\otimes a_n\xi, b_1\xi
\otimes\dots\otimes b_m\xi\rangle &=& \delta_{m,n}\langle
a_n\xi,\langle\dots,\langle a_1\xi, b_1\xi\rangle
b_2\xi\rangle,\dots b_n\xi\rangle.
\end{eqnarray*}
($a,a_j,b,b_j\in\mathfrak{B}, j=1,\dots,n$)

Note that the $\mathfrak{B}$-linear operators $A_1,
A_2:\mathcal{F}\langle\xi\rangle
\longrightarrow\mathcal{F}\langle\xi\rangle$ described by the
relations
\begin{eqnarray*}
A_1b &=&
\xi b\\
A_1(a_1\xi\otimes \dots\otimes a_n\xi b) &=&
 \xi\otimes a_1\xi\otimes \dots\otimes a_n\xi b\\
A_2b&=&0\\
A_2(a_1\xi\otimes \dots\otimes a_n\xi b) &=&
\langle\xi,a_1\xi\rangle a_2 \xi\otimes \dots\otimes a_n\xi b\\
\end{eqnarray*}
are self-adjoint to each other, in the sense that
\[
 \langle A_1\widetilde{\zeta_1},\widetilde{\zeta_2}\rangle
  =
\langle\widetilde{\zeta_1},A_2\widetilde{\zeta_2}\rangle
 \]
for any
$\widetilde{\zeta_1},\widetilde{\zeta_2}\in\mathcal{F}\langle\xi\rangle$,
therefore $S=A_1+A_2$ is selfadjoint.

Moreover, for any $a,b\in\mathfrak{B}$,
\begin{eqnarray*}
\langle 1, Sa^*bS1\rangle
 &=&
  \langle aS1, bS1\rangle\\
&=&
 \langle a(A_1+A_2)1, b(A_1+A_2)1\rangle\\
&=&
 \langle a\xi, b\xi\rangle
\end{eqnarray*}
and the conclusion follows by setting
$\varphi(p(\xi))=\langle1,p(S)1\rangle$ for all
$p\in\mathfrak{B}\langle\xi\rangle$.
\end{proof}

\begin{corollary}\label{corposit}The maps $\mu$ and $\nu$ from Theorem \ref{clthm} are
positive if and only if for any $b\in\mathfrak{B}$ one has that
$\Phi(Xb^*bX)\geq0$ and $\Psi(Xb^*bX)\geq0$.
\end{corollary}
\begin{proof}
 One implication is trivial, since, if $\nu$ and $\mu$ are positive, then
 \[
  \Psi(Xb^*bX)=\nu(Xb^*bX)=\nu((bX)^*bX)\geq0
  \]
 and
 \[
  \Phi(Xb^*bX)=\mu(Xb^*bX)=\mu((bX)^*bX)\geq0.
   \]

 Suppose now that $\Phi(Xb^*bX)\geq0$ and $\Psi(Xb^*bX)\geq0$ for all $b\in\mathfrak{B}$.
 We will use the same argument as in \cite{speicher1} and \cite{mvp}.

 Consider the set of selfadjoint symbols $\{\xi_i\}_{i\geq1}$.
  On each $\mathfrak{B}$-bimodule
 $\mathfrak{B}\xi_i\mathfrak{B}$ we have the positive $\mathfrak{B}$-sesquilinear pairings
 $\langle\cdot,\cdot\rangle_\Phi$
 and $\langle\cdot,\cdot\rangle_\Psi$ determined by
 \begin{eqnarray*}
 \langle a\xi_i,b\xi_i\rangle_\Phi
 &=&
 \Phi(Xa^*bX)\\
 \langle
 a\xi_i,b\xi_i\rangle_\Psi
 &=&
 \Psi(Xa^*bX).\\
 \end{eqnarray*}

As shown in Lemma \ref{ext1}, the above $\mathfrak{B}$-sesquilinear
pairings determine positive conditional expectations $\varphi_1,
\psi_i:\mathfrak{A}_i\longrightarrow\mathfrak{B}$, where
$\mathfrak{A}_i=\mathfrak{B}\langle\xi_i\rangle$ be the
$\ast$-algebras of polynomials in $\xi$ with coefficients from
$\mathfrak{B}$, $i\geq1$.

For $\tau:\mathfrak{B}\langle\xi\rangle\longrightarrow \mathfrak{B}$
a conditional expectation, and $\lambda\geq 0$, note with
$D_\lambda\tau$ the dilation with $\lambda$ of $\tau$, i.e.
  \[
 D_\lambda\tau(\xi b_1\xi\cdots b_n\xi)
  =
   \lambda^{n+1}\tau(\xi b_1\xi\cdots b_n\xi)
   \]
Remark that if $\tau$ is positive, then $D_\lambda\tau$ is also
positive.

 With the notations above, consider, as in Definition \ref{defncfree}, the conditionally
free product $(\mathfrak{A},
\Phi,\Psi)=\free_{i\in\mathcal{I}}(\mathfrak{A}_i,\Phi_i,\Psi_i)$.
The elements $\{\xi_i\}_{i\geq1}$ are conditionally free in
$\mathfrak{A}$, so Theorem \ref{clthm} implies that:

\begin{eqnarray*}
\mu
 &=&
  \lim_{N\longrightarrow\infty}\Phi_{\frac{\xi_1+\dots +\xi_N}{\sqrt{N}}}
=D_{\frac{1}{\sqrt{N}}}\Phi_{\xi_1+\dots +\xi_N}\\
\nu
 &=&
  \lim_{N\longrightarrow\infty}\Psi_{\frac{\xi_1+\dots +\xi_N}{\sqrt{N}}}
 =D_{\frac{1}{\sqrt{N}}}\Psi_{\xi_1+\dots +\xi_N}\\
&=&
 D_{\frac{1}{\sqrt{N}}}\left(\free_{i=1}^N\Psi_{\xi_i}\right).
\end{eqnarray*}

We have that $\free_{i=1}^N\Psi_{\xi_i}\geq0$ since it is the free
product of states (see, for example \cite{speicher1}), hence the
positivity of $\nu$.

Also, Theorem \ref{assoc} and Corollary \ref{sumposit} imply that
$\Phi_{\xi_1+\dots +\xi_N}\geq0$,  therefore $\mu\geq0$.

\end{proof}

\par\bigskip\noindent
{\bf Acknowledgment.} This research was partially supported by the
Grant 2-CEx06-11-34 of the Romanian Government. I am thankful to
Marek Bo\.{z}ejko for presenting me the basics of c-freeness and
bringing to my attention the references \cite{bls} and
\cite{mlotk1}. I thank also Hari Bercovici for his constant support
and his many advices during the work on this paper.

\bibliographystyle{amsplain}

%\noindent $\clubsuit$ Note to author:
%Proceedings articles should be formatted as in reference 1
%above, journal articles as in reference 2 above, and books as in
%reference 3 above.

\end{document}